\documentclass{article}

\usepackage[ansinew]{inputenc}
\usepackage{theorem}
\usepackage[english]{babel}
\usepackage{graphicx}
\usepackage{amsfonts}
\usepackage{amsmath}
\usepackage[dvips]{color}

\newtheorem {defn}{Definition}

\newtheorem {thm}{Theorem}

\newtheorem {lemma}{Lemma}

\newcommand{\SP}{{\operatorname{sp}}}
\newcommand{\TR}{{\operatorname{tr}}}
\newcommand{\SPEC}{{\operatorname{spec}}}
\newcommand{\diag}{\operatorname{diag}}
\newcommand{\bbmatrix}[1]{\left[ \begin{array}{cccccccccccccccccc} #1 \end{array} \right]}
\def\squarebox#1{\hbox to #1{\hfill\vbox to #1{\vfill}}}
\newcommand{\qed}{\hspace*{\fill}
\vbox{\hrule\hbox{\vrule\squarebox{.667em}\vrule}\hrule}\smallskip}

\title{On construction of solutions of evolutionary Non Linear Schr\" odinger equation}
\author{Andrey Melnikov\\Drexel university, Philadelphia, USA}

\begin{document}
\maketitle

\abstract{In this work we present an application of a theory of vessels to solution of the evolutionary Non Liner Schr\" odinger (NLS) equation.
The classes of functions for which the initial value problem is solvable relies on the existence of an analogue of the inverse scattering theory
for the usual NLS equation. This approach is similar to the classical approach of Zackarov-Shabbath for solving of evolutionary
NLS equation, but has an advantage of simpler formulas and new techniques and notions to construct solutions of the evolutionary NLS equation.
}

\tableofcontents
\section{Introduction and Background}
Solution of the Non Lineat Schr\" odinger (NLS) evolutionary equation plays a special role in the theory of PDEs and in physics (optics and water waves).
This equation can be defined as follows
\begin{equation} \label{eq:ENLS}
i y_t + y_{xx} + 2 |y|^2 y = 0, \quad y(x,0) = \beta(x),
\end{equation}
where $y=y(x,t)$ is a complex valued function of two real variables $x,t$ and $y(x,0)=\beta(x)$ is the initial condition, defined on an interval
$\mathrm I\subseteq\mathbb R$.
Notice that constant $2$ used in this work can be replaced with arbitrary number but scaling. This equation has so called integrability property,
which enables to use the inverse scattering theory in order to solve it. It is usually done using Zakharov- Shabat system
\cite{bib:ZakSha74}. There are also numerous numerical solutions of this equation, among which we can mention split-step (Fourier) method.
We are going to generalize the method of Zakharov-Shabat by introducing evolutionary NLS (regular) vessels. The setting we choose to work in
involves bounded operators only and Hilbert space techniques. It is worth noticing that one can also generalize this theory to unbounded operators
(as it was done in \cite{bib:GenVessel} for the Sturm-Liuoville differential equation).

We assume from the beginning that the initial value $\beta(x)$ is analytic on $\mathrm I\subseteq\mathbb R$ function and moreover arises from a regular
NLS vessel. It is on going research, whether every analytic on $\mathrm I$ function can be presented using a vessel and we tempt to believe
that it is the case. The problem of construction of a vessel for a given $\beta(x)$ is still open and complicated. On the other hand, we show examples
of how to construct different $\beta(x)$, which have a bounded spectrum, a spectrum on a bounded, continuous curve, or an infinite discrete spectrum.
The problem of constructing of a vessel for a given $\beta(x)$ is identical to the inverse scattering problem of the NLS equation. To emphasize this
point let us briefly discuss how the construction of the solution of \eqref{eq:ENLS} is performed starting from the \textit{spectral data} encoded
in a $2\times 2$ matrix-function $S(\lambda)$ which is realized \cite{bib:bgr} in the following manner
\[ S(\lambda) = I - B_0^* \mathbb X_0^{-1}(\lambda I - A)^{-1}B_0.
\]
Here $\mathcal H$ is an auxiliary Hilbert space, on which there are defined bounded operators $A, \mathbb X_0$ and the operator $B_0:\mathbb C^2\rightarrow\mathcal H$.
We also assume that
\[ \mathbb X_0^*=\mathbb X_0, \quad A \mathbb X_0 + \mathbb X_0 A^* + B_0B_0^* = 0.
\]
Then solving for $B(x)$ the equation \eqref{eq:DB} $B'(x) = -AB(x)\sigma_2$, for $\sigma_2=\dfrac{1}{2} \bbmatrix{1&0\\0&-1}$ and defining next
$\mathbb X(x) = \mathbb X_0 + \int_{x_0}^x B(y)\sigma_2B^*(y)dy$ we will obtain that the function
\[
\beta(x) = \bbmatrix{1&0}B^*\mathbb X^{-1}B \bbmatrix{0\\1}
\]
is analytic (Theorem \ref{thm:betaAnalytic}) on an interval, where $\mathbb X(x)$ is invertible. On the other hand, $\mathbb X(x)$ is invertible on an interval $\mathrm I$ including $x_0$,
because $\mathbb X_0$ is such. Indeed, at the definition of $\mathbb X(x)$ we add to $\mathbb X_0$ a bounded operator with bounded growth. Moreover, it turns out
(Theorem \ref{thm:Backlund}) that the function
\[ S(\lambda,x) = I - B^*(x)\mathbb X^{-1}(x)(\lambda I-A)^{-1} B(x)
\]
maps solutions $u(\lambda,x)$ of the trivial NLS equation \eqref{eq:InCC}
\[ \dfrac{\partial}{\partial x} u(\lambda,x) - (- i \lambda) \sigma_2 u(\lambda,x) = 0
\]
to solutions $y(\lambda,x) = S(\lambda,x) u(\lambda,x)$ of the NLS equation \eqref{eq:OutCC}
\[
 \dfrac{\partial}{\partial x} y(\lambda,x) - (- i \lambda) \sigma_2 (\lambda,x) + 
\bbmatrix{0&\beta(x)\\-\beta^*(x)&0} u(x,\lambda) = 0.
\]
In order to obtain a solution of the evolutionary NLS \eqref{eq:ENLS},
we evolve the operator $B$ with respect to $t$ as follows 
\[ \dfrac{\partial}{\partial t} B(x,t) = i A \dfrac{\partial}{\partial x} B(x,t) = i A (-A B(x,t) \sigma_2), \quad B(x,0) = B(x)
\]
and redefine
\[ \mathbb X(x,t) = \mathbb X(x) + \int_0^t [i A B(x,s) \sigma_2 B^*(x,s) - i B(x,s) \sigma_2 B^*(x,s) A^*] ds.
\]
One of the most interesting results of this paper is Theorem \ref{thm:EbeatEqn}, where we prove that the new
\[ \beta(x,t) =  \bbmatrix{1&0}B^*(x,t)\mathbb X^{-1}(x,t)B(x,t) \bbmatrix{0\\1}
\]
satisfies \eqref{eq:ENLS} and coincides with $\beta(x)$ for $t=0$, proving the existence of solutions for \eqref{eq:ENLS} with initial value.

The formulas presented here enable also to explicitly perform this construction for some basic and important cases. We show how to construct
a vessel, for which the spectrum of $A$ lies on a curve $\Gamma$ (Section \ref{sec:VesselFromGamma} for $\beta(x)$ and \ref{sec:EVesselFromGamma}
for $\beta(x,t)$), with a discrete set (Section \ref{sec:VesselFromD} for $\beta(x)$ and Section \ref{sec:EVesselFromD} for $\beta(x,t)$).
We also discuss the general construction in Sections \ref{sec:VesselFromS0} and \ref{sec:EVesselFromS0}, Finally we present constructions of the Solitons
in Section \ref{sec:Solitons}.

\section{Definition of a regular NLS vessel and its properties}
We define first parameters, which will be frequently used in the sequel. Other choice of these parameters generates solutions of the Sturm-Liuoville differential equation
and the Kortweg-de-Vries equation (see \cite{bib:GenVessel, bib:KdVVessels} for details).
\begin{defn} NLS vessel parameters are defined as follows:
\[ \sigma_1 = \bbmatrix{1&0\\0&1},\quad
\sigma_2 = \dfrac{1}{2} \bbmatrix{1&0\\0&-1},
\quad \gamma =\bbmatrix{0 & 0 \\0 & 0},
\]
\end{defn}
\begin{defn} 
An NLS regular vessel is a collection of operators and spaces:
\begin{equation} \label{eq:NLS}
\mathfrak{V}_{NLS} = (A, B(x), \mathbb X(x); \sigma_1, 
\sigma_2, \gamma, \gamma_*(x);
\mathcal{H},\mathbb C^2,\mathrm I),
\end{equation}
where $\sigma_1, \sigma_2, \gamma$ are the NLS vessel parameters, $\mathrm I$ is a closed interval.
The operators $A, \mathbb X(x):\mathcal H\rightarrow\mathcal H$, $B(x):\mathbb C^2\rightarrow\mathcal H$ are bounded operators for all $x\in\mathrm I$.
The operator $\mathbb X(x)$ is invertible for all $x\in\mathrm I$. The operators are subject to the
following vessel conditions:
\begin{align}
\label{eq:DB} 0  =  \frac{d}{dx} (B(x)) + A B(x) \sigma_2, \\
\label{eq:Lyapunov} A \mathbb X(x) + \mathbb X(x) A^* + B(x) B(x)^* = 0, \\
\label{eq:DX} \frac{d}{dx} \mathbb X(x)  =  B(x) \sigma_2 B(x)^*, \\
\label{eq:Linkage}
\gamma_*(x)  = \sigma_2 B(x)^* \mathbb X^{-1}(x) B(x) 
- B(x)^* \mathbb X^{-1}(x) B(x) \sigma_2, \\
\mathbb X^*(x)=\mathbb X(x),
\end{align}
\end{defn}
For each vessel NLS vessel there are exist three notions, which play a significant role in research.
\begin{defn}Suppose that we are given an NLS vessel $\mathfrak{V}_{NLS}$. Then
its \textbf{transfer} function $S(\lambda,x)$, the \textbf{tau}-functions $\tau(x)$ and the \textbf{beta}-function $\beta(x)$ are defined as follows:
\begin{eqnarray} 
\label{def:S} S(\lambda, x) = I - B(x)^* \mathbb X^{-1}(x) (\lambda I - A)^{-1} B(x), \\
\label{eq:DefTau} \tau(x) = \det\big(\mathbb X^{-1}(x_0)\mathbb X(x)\big),\\
\label{eq:DefBeta} \beta(x) = \bbmatrix{1&0}\gamma_*(x)\bbmatrix{0\\1} = \bbmatrix{1&0}B^*(x)\mathbb X^{-1}(x)B(x) \bbmatrix{0\\1}.
\end{eqnarray}
\end{defn}
The definition of $\beta(x)$ may be considered as excessive, because actually the matrix-function $\gamma_*(x)$ turns to be
\[ \gamma_*(x) = \bbmatrix{0&\beta(x)\\-\beta^*(x)&0},
\]
using the self-adjointenss of $\mathbb X(x)$. Still, we will use these two notions extensively, so we have defined both of them.
Notice that $S(\lambda,x)$ is a $2\times 2$ matrix-function, whose poles and singularities with respect to $\lambda$ are determined by the operator$A$ only.
Since all the involved operators, appearing at the definition of $S(\lambda,x)$ are
bounded, we can see that $S(\lambda,x)$ is analytic in $\lambda$ for all $x\in\mathrm I$ with value $I$ (=identity $2\times 2$ matrix) there. As a result,
we can consider its Taylor series
\[ S(\lambda,x) = I - B^*(x) \mathbb X^{-1}(x) (\lambda I - A)^{-1} B(x)\sigma_1 =
I - \sum\limits_{i=0}^\infty \dfrac{H_i(x)}{\lambda^{i+1}} \sigma_1,
\]
which is convergent at least for $\lambda > \| A \|$. Thus we define its (Markov) moments as follows
\begin{defn}
$n$-th moment of the vessel $\mathfrak{V}_{NLS}$ is  
\[ H_n(x) = B^*(x) \mathbb X^{-1}(x) A^n B(x).
\]
\end{defn}
We will present in the next section basic properties of an NLS vessel $\mathfrak{V}_{NLS}$ by exploring all the objects (the transfer, the tau, the beta functions
and the moments). We will also see that there is a standard technique for construction of such vessels.

\subsection{The transfer and the tau function of an NLS vessel}
The main reason to consider NLS vessel is the next Theorem \ref{thm:Backlund}, whose proof we
present in full details, also it appeared in different settings in \cite{bib:defVess, bib:Vortices, bib:MyThesis,bib:MelVin1, bib:MelVinC, bib:SLVessels}.
It turns out that we can see an NLS vessel as a B\" acklund transformation of the trivial NLS equation to a more complicated one. More precisely, the transfer
function $S(\lambda,x)$ of such a vessel maps solutions of the so called \textit{input LDE} with the spectral parameter $\lambda$
\begin{equation} \label{eq:InCC}
[\lambda \sigma_2 - \sigma_1 \dfrac{\partial}{\partial x} + \gamma] u(\lambda,x) = 0
\end{equation}
to solutions of the \textit{output LDE} with the same spectral parameter
\begin{equation} \label{eq:OutCC}
[\lambda \sigma_2 - \sigma_1 \dfrac{\partial}{\partial x} + \gamma_*] y(\lambda,x) = 0
\end{equation}
As a result the following differential equation holds
\begin{equation} \label{eq:DS}
S'_x = \sigma_1^{-1}(\sigma_2 \lambda + \gamma_*) S - S \sigma_1^{-1}(\sigma_2 \lambda + \gamma).
\end{equation}

\begin{thm}[Vessel as a B\" acklund transformation] \label{thm:Backlund} Suppose that $\mathfrak{V}_{NLS}$ is an NLS vessel, defined in \eqref{eq:NLS} and suppose
that $u(x,\lambda)=\bbmatrix{u_1(\lambda,x)\\u_2(\lambda,x)}$ be a solution of the input LDE \eqref{eq:InCC}
\[ \lambda \sigma_2 u(\lambda, x) -
\sigma_1 \frac{\partial}{\partial x}u(\lambda,x) +
\gamma u(\lambda,x) = 0,
\]
then $y(\lambda,x)=S(\lambda,x)u(\lambda)$ is a solution of the output LDE \eqref{eq:OutCC}: 
\[
\lambda \sigma_2 y(\lambda, x) - \sigma_1 \frac{\partial}{\partial x}y(\lambda,x) +
\gamma_*(x) y(\lambda,x) = 0.
\]
\end{thm}
\noindent\textbf{Proof:} we plug in the expression
\[ y(\lambda,x)=S(\lambda,x)u(\lambda) = (I - B^*(x)\mathbb X^{-1}(x)(\lambda I - A)^{-1} B(x)\sigma_1) u(\lambda,x)
\]
into \eqref{eq:OutCC} for all $\lambda\not\in\SPEC(A)$. Let us denote
$G(\lambda,x)=B^*(x)\mathbb X^{-1}(x)(\lambda I - A)^{-1} B(x)$, then for all $x\in\mathrm I$ and $\lambda\not\in\SPEC(A)$
\[ \begin{array}{lllllll}
\lambda \sigma_2 y(\lambda, x) - \sigma_1 \frac{\partial}{\partial x}y(\lambda,x) +
\gamma_*(x) y(\lambda,x) = \\
\quad = \lambda \sigma_2 [(I - G(\lambda,x)\sigma_1) u(\lambda,x)] - \sigma_1 \frac{\partial}{\partial x} [(I - G(\lambda,x)\sigma_1) u(\lambda,x)] + \\
\quad \quad \quad + \gamma_*(x) (I - G(\lambda,x)\sigma_1 u) (\lambda,x) = \\
\quad = \lambda \sigma_2 u(\lambda,x) - \lambda \sigma_2 G(\lambda,x) \sigma_1 u(\lambda,x) - \\
\quad \quad - \sigma_1 \frac{\partial}{\partial x} u(\lambda,x) + \sigma_1 \frac{\partial}{\partial x}[ G(\lambda,x)] \sigma_1 u(\lambda,x) -
 \sigma_1 G(\lambda,x)\sigma_1 \frac{\partial}{\partial x} u(\lambda,x) + \\
\quad \quad \quad + \gamma_*(x) u(\lambda,x) - \gamma_*(x)G(\lambda,x) \sigma_1 u(\lambda,x) = \\
\quad = \text{using \eqref{eq:InCC}} = \\
\quad = - \lambda \sigma_2 G(\lambda,x) \sigma_1 u(\lambda,x) + \sigma_1 \frac{\partial}{\partial x}[ G(\lambda,x)] \sigma_1 u(\lambda,x) - 
 \sigma_1 G(\lambda,x) \sigma_1 \frac{\partial}{\partial x} u(\lambda,x) + \\
\quad \quad \quad + (\gamma_*(x)-\gamma) u(\lambda,x) - \gamma_*(x) G(\lambda,x) \sigma_1 u(\lambda,x).
\end{array} \]
Let us differentiate the expression $\sigma_1 G(\lambda,x) \sigma_1$ using formulas \eqref{eq:DB} and \eqref{eq:DX}:
\[ \begin{array}{lllllll}
\sigma_1 \frac{\partial}{\partial x} G(\lambda,x) \sigma_1 = \\
= \frac{\partial}{\partial x}[\sigma_1 B^*(x)\mathbb X^{-1}(x)(\lambda I - A)^{-1} B(x)\sigma_1] = \\
= (-\sigma_2 B^*(x) A^* - \gamma B^*(x)) \mathbb X^{-1}(x)(\lambda I - A)^{-1} B(x)\sigma_1 - \sigma_1 B^*(x)\mathbb X^{-1}(x) B(x) \sigma_2 G(\lambda,x) \sigma_1 + \\
\quad \quad + \sigma_1 B^*(x)\mathbb X^{-1}(x)(\lambda I - A)^{-1} (-A B(x)\sigma_2 - B(x)\gamma) = \\
= -\sigma_2 B^*(x) A^* \mathbb X^{-1}(x)(\lambda I - A)^{-1} B(x)\sigma_1 -
\gamma G(\lambda,x) \sigma_1 - \sigma_1 B^*(x)\mathbb X^{-1}(x) B(x) \sigma_2 G(\lambda,x) \sigma_1 - \\
\quad \quad - \sigma_1 B^*(x)\mathbb X^{-1}(x)(\lambda I - A)^{-1} A B(x)\sigma_2 - \sigma_1 G(\lambda,x) \gamma  = \\
= \text{using \eqref{eq:Lyapunov} for $ A^* \mathbb X^{-1}(x)$ } = \\
= -\sigma_2 B^*(x) (-\mathbb X^{-1}(x) A - \mathbb X^{-1}(x)B(x) \sigma_1 B^*(x)\mathbb X^{-1}(x)) (\lambda I - A)^{-1} B(x)\sigma_1 - \\
\quad - \gamma G(\lambda,x) \sigma_1 - \sigma_1 B^*(x)\mathbb X^{-1}(x) B(x) \sigma_2 G(\lambda,x) \sigma_1 - \\
\quad \quad - \sigma_1 B^*(x)\mathbb X^{-1}(x)(\lambda I - A)^{-1} A B(x)\sigma_2 - \sigma_1 G(\lambda,x) \gamma  = \\
= \sigma_2 B^*(x) \mathbb X^{-1}(x) (A\pm \lambda I) (\lambda I - A)^{-1} B(x)\sigma_1 +  \sigma_2 B^*(x) \mathbb X^{-1}(x) B(x) \sigma_1 G(\lambda,x) \sigma_1 -\\
\quad - \gamma G(\lambda,x) \sigma_1 - \sigma_1 B^*(x)\mathbb X^{-1}(x) B(x) \sigma_2 G(\lambda,x) \sigma_1 - \\
\quad \quad - \sigma_1 B^*(x)\mathbb X^{-1}(x)(\lambda I - A)^{-1} (A\pm\lambda I) B(x)\sigma_2 - \sigma_1 G(\lambda,x) \gamma  = \\
= - \sigma_2 B^*(x) \mathbb X^{-1}(x) B(x)\sigma_1 +  \lambda \sigma_2 G(\lambda,x) \sigma_1 + 
+  \sigma_2 B^*(x) \mathbb X^{-1}(x) B(x) \sigma_1 G(\lambda,x) \sigma_1 -\\
\quad - \gamma G(\lambda,x) \sigma_1 - \sigma_1 B^*(x)\mathbb X^{-1}(x) B(x) \sigma_2 G(\lambda,x) \sigma_1 + \\
\quad \quad + \sigma_1 B^*(x)\mathbb X^{-1}(x) B(x)\sigma_2 
- \lambda \sigma_1G(\lambda,x) \sigma_2  - \sigma_1 G(\lambda,x) \gamma.
\end{array} \]
Plugging this expression back into the equation \eqref{eq:OutCC} developed earlier and performing some obvious cancellations,
we will obtain
\[ \begin{array}{lllllll}
\lambda \sigma_2 y(\lambda, x) - \sigma_1 \frac{\partial}{\partial x}y(\lambda,x) +
\gamma_*(x) y(\lambda,x) = \\
= - \lambda \sigma_2 G(\lambda,x) \sigma_1 u(\lambda,x) + \sigma_1 \frac{\partial}{\partial x}[ G(\lambda,x)] \sigma_1 u(\lambda,x) - 
 \sigma_1 G(\lambda,x) \sigma_1 \frac{\partial}{\partial x} u(\lambda,x) + \\
\quad \quad \quad + (\gamma_*(x)-\gamma) u(\lambda,x) - \gamma_*(x) G(\lambda,x) \sigma_1 u(\lambda,x) = \\
= (-\sigma_2B^*(x)\mathbb X^{-1}(x)B(x)\sigma_1+\sigma_1B^*(x)\mathbb X^{-1}(x)B(x)\sigma_1+\gamma_*(x)-\gamma)  u(\lambda,x) + \\
\quad  + (\sigma_2B^*(x)\mathbb X^{-1}(x)B(x)\sigma_1-\sigma_1B^*(x)\mathbb X^{-1}(x)B(x)\sigma_1+\gamma_*(x)-\gamma)  G(\lambda,x)\sigma_1 u(\lambda,x) + \\
\quad \quad + \sigma_1  G(\lambda,x) ( \sigma_1 \frac{\partial}{\partial x}  u(\lambda,x) - \lambda \sigma_2  u(\lambda,x) + \gamma u(\lambda,x)) = 0,
\end{array} \]
using the linkage condition \eqref{eq:Linkage} and the differential equation \eqref{eq:InCC}. \qed

Let us present next the significance (and well definedness) of the tau-function $\tau(x)$, defined in \eqref{eq:DefTau}. 
Using vessel condition \eqref{eq:DX} $\mathbb X(x)$ has the formula
\[ \mathbb X(x) = \mathbb X(x_0) + \int\limits_0^x B^*(y) \sigma_2 B(y) dy,
\]
and as a result
\[ \mathbb X^{-1}(x_0) \mathbb X(x) = I + \mathbb X^{-1}(x_0) \int\limits_{x_0}^x B^*(y) \sigma_2 B(y) dy.
\]
Since $\sigma_2$ has rank 2, this expression is of the form $I + T$, for a trace class operator $T$ and since 
$\mathbb X_0$ is an invertible operator, there exists a non trivial interval (of length at least $\dfrac{1}{\|\mathbb X_0^{-1}\|}$) on which $\mathbb X(x)$ and $\tau(x)$ are defined. Recall \cite{bib:GKintro} that a function $F(x)$ from (a, b) into the group G (the set of bounded invertible operators on H of the form I + T, for
a trace-class operator $T$) is said to be differentiable if $F(x) -I$ is \textit{differentiable} as a map into the trace-class operators. In our case,
\[ \dfrac{d}{dx} (\mathbb X^{-1}(x_0)\mathbb X(x)) = 
\mathbb X^{-1}(x_0) \dfrac{d}{dx} \mathbb X(x) =
\mathbb X^{-1}(x_0) B(x)\sigma_2 B^*(x)
\]
exists in trace-class norm. Israel Gohberg and Mark Krein \cite[formula 1.14 on p. 163]{bib:GKintro}
proved that if $\mathbb X^{-1}(x_0)\mathbb X(x)$ is a differentiable function
into G, then $\tau(x) = \SP (\mathbb X^{-1}(x_0)\mathbb X(x))$
\footnote{$\SP$ - stands for the trace in the infinite dimensional space.} is a differentiable map into $\mathbb C^*$ with
\begin{multline} \label{eq:tauDet}
\dfrac{\tau'}{\tau}  = \SP (\big(\mathbb X^{-1}(x_0)\mathbb X(x)\big)^{-1} 
\dfrac{d}{dx} \big(\mathbb X^{-1}(x_0)\mathbb X(x)\big)) = \SP (\mathbb X(x)' \mathbb X^{-1}(x)) = \\
= \SP (B(x)\sigma_2 B^*(x) \mathbb X^{-1}(x)) =
\TR (\sigma_2 B^*(x) \mathbb X^{-1}(x)B(x)).
\end{multline}
A most important question, related to this theory is for what classes of $\beta(x)$ are obtained. This question is answered in the next Theorem.
\begin{thm} \label{thm:betaAnalytic}
Suppose that $\mathfrak{V}_{NLS}$ is an NLS vessel, defined in \eqref{eq:NLS}. Then the function $\beta(x)$ is analytic on the interval $\mathrm I$.
\end{thm}
\noindent\textbf{Proof:} Notice that from the formula \eqref{eq:DB} it follows the operator $B(x)$ is analytic in $x$. Since $\mathbb X(x)$ is invertible on $\mathrm I$,
the operator $\mathbb X^{-1}(x)$ is also analytic in $x$ using the formula
\[ \dfrac{d}{dx} \mathbb X^{-1}(x) = - \mathbb X^{-1}(x) B(x)\sigma_2B^*(x) \mathbb X^{-1}(x)
\]
Thus $\beta(x)$, defined by \eqref{eq:DefBeta} is analytic on $\mathrm I$.
\qed

\begin{thm}[Permanency conditions] \label{thm:Permanency}
Suppose that we are given an NLS regular vessel $\mathfrak{V}_{NLS}$, then
\begin{enumerate}
	\item if the Lyapunov equation \eqref{eq:Lyapunov} holds for a fixed $x_0\in\mathrm I$, then it holds for all $x\in\mathrm I$,
	\item if $S^*(-\bar\lambda,x_0) \sigma_1S(\lambda,x_0) = \sigma_1$ then 
	\begin{equation} \label{eq:Symmetry}
	S^*(-\bar\lambda,x) \sigma_1 S(\lambda,x) = \sigma_1 
	\end{equation}
	for all $x\in\mathrm I$,
	\item $\det S(\lambda,x) = \det S(\lambda,x_0)$ for all $x_0, x\in\mathrm I$ and all points of $\lambda$-analyticity of $S(\lambda,x)$.
\end{enumerate}
\end{thm}
\noindent\textbf{Proof:} Differentiating Lyapunov equation and using the vessel conditions \eqref{eq:DB}, \eqref{eq:DX} we will obtain that
\begin{multline}
 A \dfrac{d}{dx} \mathbb X(x) + \dfrac{d}{dx} \mathbb X(x) A^* + \dfrac{d}{dx} [B(x)] B^*(x) + B(x)  \dfrac{d}{dx} B^*(x) = \\
A B(x)\sigma_2B^*(x) +  B(x)\sigma_2B^*(x)A^* + (-AB(x)\sigma_2)B^*(x) + B(x) (-\sigma_2 B^*A^*) = 0,
\end{multline}
from where the permanency of the Lyapunov equation follows. Similarly, differentiating $S^*(-\bar\lambda,x) \sigma_1 S(\lambda,x)$ we will obtain zero
and the permanency of \eqref{eq:Symmetry} follows. For the last statement, using \eqref{eq:DS} we calculate for $\lambda\not\in\SPEC(A)$
\[ \begin{array}{llll}
\dfrac{\dfrac{\partial}{\partial x} \det S(\lambda,x)}{\det S(\lambda,x)} =
\TR (S^{-1}(\lambda,x)\dfrac{\partial}{\partial x} S(\lambda,x)) = \\
= \TR(S^{-1}(\lambda,x) [(\sigma_2 \lambda + \gamma_*(x)) S(\lambda,x)
-  S(\lambda,x)(\sigma_2 \lambda + \gamma)] ) = \\
= \TR(\sigma_2 \lambda + \gamma_*(x) 
- (\sigma_2 \lambda + \gamma)) = \TR(\gamma_*(x) - \gamma(x)) = \\
= \TR(B^*(x)\mathbb X^{-1}(x) B(x) \sigma_2 -
\sigma_2 B^*(x)\mathbb X^{-1}(x) B(x)) = 0,
\end{array} \]
\qed

We will not be using the last property in this statement, but we find it interesting by itself.

\subsection{Moments}
The following properties of the moments $H_(x)$ of an NLS vessel are immediate from their definition as the coefficients of
$\dfrac{1}{\lambda^{n+1}}$ at the Taylor series of $S(\lambda,x)$.
\begin{thm} Let $\mathfrak{V}_{NLS}$ be an NLS vessel. Then its moments satisfy the following equations
\begin{eqnarray} 
\label{eq:HiHi+1} \sigma_1^{-1} \sigma_2 H_{n+1} - H_{n+1} \sigma_2 \sigma_1^{-1} =
 \frac{d}{dx} H_n  - \sigma_1^{-1} \gamma_* H_n  + H_n \gamma \sigma_1^{-1}, \\
\label{eq:Moments} H_{n+1} + (-1)^n H_{n+1}^*  = \sum_{j=0}^{n} (-1)^{j+1} H_{n-j} \sigma_1 H_j^*.
\end{eqnarray}
\end{thm}
\noindent\textbf{Proof:} Plugging the Taylor expansion formula into \eqref{eq:DS} and equating the coefficients of $\dfrac{1}{\lambda^{n=1}}$
we will obtain the first formula \eqref{eq:HiHi+1}. The second formula \eqref{eq:Moments} is obtained in the same manner from \eqref{eq:Moments}.
\qed

It turns out that using only the differential equations \eqref{eq:HiHi+1} one can create a recursive formula for the entries
$H_n = \bbmatrix{H_n^{11}&H_n^{12}\\H_n^{21}&H_n^{22}}$ as follows
\begin{equation} \label{eq:NLSmoments} 
\left\{ \begin{array}{llll}
H_{n}^{21} = - \dfrac{d}{dx}\big( H_{n-1}^{21}\big) - \beta^* H_{n-1}^{11}, \\
H_{n}^{12} = \dfrac{d}{dx}\big( H_{n-1}^{12}\big) - \beta H_{n-1}^{22}, \\
\dfrac{d}{dx}\big( H_n^{11}\big) = \beta H_n^{21}, \\
\dfrac{d}{dx}\big( H_n^{22}\big) = -\beta^* H_n^{12}.
\end{array} \right. \end{equation}
while the first moment $H_0$ is found from the linkage condition \eqref{eq:Linkage}
\[ H_0 = \bbmatrix{ H_0^{11} & \beta \\ \beta^* & H_0^{22}},
\]
and the entries $H_0^{11}, H_0^{22}$ are found using two last equations of \eqref{eq:NLSmoments}:
\[ \begin{array}{ll}
\dfrac{d}{dx}\big( H_0^{11}\big) = \beta H_0^{21} = \beta \beta^* , \\
\dfrac{d}{dx}\big( H_0^{22}\big) = -\beta^* H_0^{12} = -\beta^* \beta .
\end{array} \]
As a result, we obtain that the formula for the tau-function \eqref{eq:tauDet} becomes
\[ \dfrac{\tau'}{\tau}  = \TR (\sigma_2 B^*(x) \mathbb X^{-1}(x)B(x)) =
\TR (\sigma_2 H_0) = \dfrac{1}{2}(H_0^{11}-H_0^{22}) = 
\dfrac{1}{2}(H_0^{11}(0)-H_0^{22}(0)) + \int_{x_0}^x |\beta(y)|^2dy.
\]
For example, if we choose the initial parameters $H_0^{11}(0),H_0^{22}(0)$ to be equal, we will obtain that
under normalization $\tau(x_0)=1$:
\begin{equation} \tau(x) = e^{\int_{x_0}^x \int_{x_0}^t |\beta(y)|^2dy dt}
\end{equation}

\section{Examples of constructions of regular NLS vessels}
\subsection{\label{sec:VesselFromS0}Construction of an NLS vessel from a realized function}
Construction of an NLS vessel from the scattering data (initial condition $S(x_0,\lambda)$) can be performed as follows.
Suppose that the function $S(\lambda,x_0)$ is realized \cite{bib:bgr} as follows
\[ S(x_0,\lambda) = I - B_0 \mathbb X_0^{-1}(\lambda I -A)^{-1} B_0 \sigma_1,
\]
satisfying additionally  $A \mathbb X_0 + \mathbb X_0 A^* + B_0 B_0^* = 0$ and $\mathbb X_0^* = \mathbb X_0$. These two conditions
are required by the permanency conditions (Theorem \ref{thm:Permanency}) and will hold for all $x$ by the construction.
then define $B(x)$ as the unique solution of \eqref{eq:DB} with initial value $B_0$:
\begin{equation} \label{eq:DefBFromS0}
0  =  \frac{d}{dx} (B(x)) + A B(x) \sigma_2, \quad B(0)=B_0.
\end{equation}
Then define
\begin{equation} \label{eq:DefXFromS0}
 \mathbb X(x) = \mathbb X_0 + \int_{x_0}^x B(y)\sigma_2B^*(y)dy,
\end{equation}
and define $\gamma_*(x)$ (and hence $\beta(x)$) using \eqref{eq:Linkage}. It is a straigtforward to check
that all vessel conditions hold:
\begin{thm} Suppose that we are given of a function
\[ S(x_0,\lambda) = I - B_0 \mathbb X_0^{-1}(\lambda I -A)^{-1} B_0
\]
where using an auxiliary Hilbert space $\mathcal H$ the operators act as follows: 
\[ A,\mathbb X_0:\mathcal H\rightarrow\mathcal H, \quad
B_0:\mathbb C^2\rightarrow\mathcal H. \]
and satisfy $A \mathbb X_0 + \mathbb X_0 A^* + B_0 B_0^* = 0$ and $\mathbb X_0^* = \mathbb X_0$.
Define a collection \eqref{eq:NLS} using the formulas \eqref{eq:DefBFromS0}, \eqref{eq:DefXFromS0},
\eqref{eq:Linkage}
\[ 
\mathfrak{V}_{NLS} = (A, B(x), \mathbb X(x); \sigma_1, 
\sigma_2, \gamma, \gamma_*(x);
\mathcal{H},\mathbb C^2),
\]
then this collection is an NLS regular vessel existing on an interval $\mathrm I$ on which the operator $\mathbb X(x)$
is invertible.
\end{thm}
We show two special examples, arising from this construction for special cases of the choice of
$\mathcal H$.

\subsection{\label{sec:VesselFromGamma}Construction of a regular NLS vessel with the spectrum on a curve $\Gamma$}
Let us fix a bounded continuous curve $\Gamma = \{ \mu(t)\mid t\in [a,b]\}$ (i.e. $\mu(t)$ is continuous)
and define $\mathcal H = L^2(\Gamma) = \{ f(\mu) \mid \int_a^b |f(\mu(t))|^2dt < \infty\}$. We suppose without loss of generality that $x_0=0$ and we construct a vessel,
existing on an interval $\mathrm I$ including zero.

Let $A = 2\mu$, i.e. it is a bounded operator acting within $L^2(\Gamma)$. Solving differential equation
\eqref{eq:DB} we find that
\begin{equation} \label{eq:DefBonGamma}
 B(x) = \bbmatrix{e^{-\mu x} b_1(\mu) & e^{\mu x} b_2(\mu)}, \quad b_1(\mu), b_2(\mu)\in\mathcal H.
\end{equation}
Notice that the adjoint $B^*(x):\mathcal H\rightarrow\mathbb C^2$ is defined as follows
\[ B^*(x) f(\mu) = \int_a^b \bbmatrix{e^{-\bar\mu(t) x} \bar b_1(\mu(t)) \\ e^{\bar\mu(t) x} \bar b_2(\mu(t)) }f(\mu(t))dt.
\]
It is a well-defined operator, because by the Cauchy-Schwartz inequality the integrals are finite.

It turns out that the operator $\mathbb X(x)$ can be also explicitly defined as follows for any $f\in\mathcal H$. 
Notice that $\mathbb X(x) f$ is a new function at $\mathcal H$, for which we present its value at the point $\mu=\mu(t)\in\Gamma$:
\[ (\mathbb X(x) f) (\mu) = - \dfrac{1}{2} \int_a^b [\dfrac{b_1(\mu) \bar b_1(\mu(s)) e^{-(\mu+\bar\mu(s))x}}{\mu + \bar\mu(s)} +
\dfrac{b_2(\mu) \bar b_2(\mu(s)) e^{(\mu+\bar\mu(s))x}}{\mu + \bar\mu(s)}] f(\mu(s)) ds.
\]
It almost immediate that such an operator satisfies the conditions \eqref{eq:Lyapunov}, \eqref{eq:DX}.
For example for the Lyapunov equation, the expression $(A \mathbb X(x) f) (\mu) + (\mathbb X(x)A^* f)(\mu)$ is obtained by multiplying
the expression under the integral in $(\mathbb X(x)f)(\mu)$ by $2(\mu+\bar\mu(s))$. When it is canceled with the denominator, we obtain 
$(-B(x)\sigma_2B^*(x)f) (\mu)$. Similarly, differentiating $(\mathbb X(x) f) (\mu)$ we cancel the denominator and switch the sign of the first term, so
that \eqref{eq:DX} holds. there is only one problem, arising from the zero of the denominator $\mu(t) + \bar\mu(s)=0$. It can be overcome by requiring that
either $\Gamma \cap (-\Gamma^*) = \emptyset$, or $b_1,b_2$ are H\" older functions and $b_1(\mu(t))\bar b_1(\mu(s)) + b_2(\mu(t))\bar b_2(\mu(s)) = 0$, whenever $\mu(t) + \bar\mu(s)=0$
(so that the zero of the denominator is canceled by the numerator).

We will make the following assumption on the curve, to simplify arguments: 
\begin{defn} Suppose that $\Gamma \cap (-\Gamma^*) = \emptyset$, where $-\Gamma^* = \{-\bar\mu(t) \mid t\in[a,b]\}$.
\end{defn}
Then investigating the formula for $\mathbb X(x)$ we find that for each $f\in\mathcal H$
\begin{equation} \label{eq:DefXonGamma}
 (\mathbb X(x) f) (\mu) =  - \dfrac{b_1(\mu) e^{-\mu x}}{2} \int_a^b \dfrac{\bar b_1(\mu(s)) e^{-\bar\mu(s)x}}{\mu + \bar\mu(s)} f(\mu(s)) ds -
\dfrac{b_2(\mu) e^{\mu x}}{2} \int_a^b \dfrac{\bar b_2(\mu(s)) e^{\bar\mu(s) x}}{\mu + \bar\mu(s)} f(\mu(s)) ds
\end{equation}
Notice that this expression can be presented as
\[ \mathbb X(x) f = b_1 e^{-\mu x} c_1(\mu) + b_2 e^{\mu x} c_2(\mu),
\]
where the functions 
\[ c_1(\mu) =  - \dfrac{1}{2} \int_a^b \dfrac{\bar b_1(\mu(s)) e^{-\bar\mu(s)x}}{\mu + \bar\mu(s)} f(\mu(s)) ds, \quad
c_2(\mu) = - \dfrac{1}{2} \int_a^b \dfrac{\bar b_2(\mu(s)) e^{\bar\mu(s) x}}{\mu + \bar\mu(s)} f(\mu(s)) ds
\]
are analytic functions of $\mu$ in $\mathbb C \backslash (-\Gamma^*)$. Notice also that $\Gamma\subseteq \mathbb C \backslash (-\Gamma^*)$
by the assumption on $\Gamma$. Thus the functions $c_1(\mu), c_2(\mu)$
can have only isolated zeros on $\Gamma$ or to be identically zero.

To simplify a proof that $\mathbb X(x)$ is an invertible operator, let us assume that $b_1=b_2$ and that they are an analytic function of $\mu$.
Then the equality $\mathbb X(x)f = 0$ is equivalent to
\[ e^{-\mu x} \int_a^b \dfrac{\bar b_1(\mu(s)) e^{-\bar\mu(s)x}}{\mu + \bar\mu(s)} f(\mu(s)) ds +
e^{\mu x} \int_a^b \dfrac{\bar b_1(\mu(s)) e^{\bar\mu(s) x}}{\mu + \bar\mu(s)} f(\mu(s)) ds = 0
\]
Taking here $\mu$ around infinity so that $e^{-\mu x}$ is big in absolute value and $e^{\mu x}$ is small we obtain that 
$\int_a^b \dfrac{\bar b_1(\mu(s)) e^{-\bar\mu(s)x}}{\mu + \bar\mu(s)} f(\mu(s)) ds = 0$. Conversely, taking  $e^{\mu x}$ big and $e^{-\mu x}$ small the second expression
$\int_a^b \dfrac{\bar b_1(\mu(s)) e^{\bar\mu(s) x}}{\mu + \bar\mu(s)} f(\mu(s)) ds = 0$ must vanish too. 
Finally, considering values of $\mu\rightarrow\infty$, we can develop into Taylor series 
\[ \dfrac{1}{\mu+\mu(s)} = \dfrac{1}{\mu} \sum (-1)^n \dfrac{\mu^n(s)}{\mu^n}
\]
and plugging it back, the moments of the function $b_1(\mu(s)) e^{\pm\mu(s)x} f(\mu(s))$ must be zero for all $s$, or that the function $b_1(\mu(s)) e^{\pm\mu(s)x} f(\mu(s))$ is
orthogonal to a dense subset $\{ \mu^n(s) \}$ of $L^2(\Gamma)$. Thus 
$f(\mu(s))$ is identically zero and $\mathbb X(x)$ is invertible for all $x$:
\begin{lemma} If $\Gamma$ is a bounded continuous curve, satisfying $\Gamma \cap (-\Gamma^*) = \emptyset$, and
the operator $\mathbb X(x)$ defined in \eqref{eq:DefXonGamma} using analytic functions $b_1(\mu)=b_2(\mu)$, then $\mathbb X(x)$ is bounded and invertible 
for all $x\in\mathbb R$.
\end{lemma}
Finally, we obtain
\begin{thm} suppose that $\Gamma$ is a bounded continuous curve, satisfying $\Gamma \cap (-\Gamma^*) = \emptyset$, Define a collection \eqref{eq:NLS}
\[ 
\mathfrak{V}_{NLS} = (A, B(x), \mathbb X(x); \sigma_1, 
\sigma_2, \gamma, \gamma_*(x);
\mathcal{H},\mathbb C^2),
\]
where $\mathcal H = L^2(\Gamma)$, $A = 2\mu$, $\mathbb X(x)$ and $B(x)$ are defined by \eqref{eq:DefXonGamma} and by \eqref{eq:DefBonGamma} for
$|b_1(\mu)|^2 + |b_2(\mu)|^2\neq 0$. Then the collection $\mathfrak{V}_{NLS}$ is an NLS regular vessel
existing on $\mathrm I = \mathbb R$.
\end{thm}
\noindent\textbf{Proof:} By the construction the operators are well defined and satisfy the vessel condition. Since the operator $\mathbb X(x)$ is invertible
for all $x\in\mathbb R$, we obtain that $\mathrm I = \mathbb R$. \qed

\subsection{\label{sec:VesselFromD}Construction of a regular NLS vessel with a discrete spectrum}
In this section we want to show how to construct a vessel, whose spectrum is a given set of numbers
$D = \{ 2 \mu_n \}$. We define $\mathcal H = \ell^2$, which the set of infinite sequences, summable in absolute value.
Now we can imitate the construction of the vessel on a curve $\Gamma$ using discretization as follows. we define first the operator
$A = \diag(2\mu_n)$ and for this operator to be bounded, we have to demand that the sequence $D$ is bounded from below and from above, namely
\begin{defn} The sequence $D$ is called \textbf{bounded} if there exist $M>0$ such that $|\mu_n|<M$ for $\mu_n\in D$. It is called
\textbf{separated from zero} if there exist $m>0$ such that $0<m<|\mu_n|$ for all $\mu_n\in D$.
\end{defn}

In the next definition we think of $\mathbb X(x)$ as an infinite matrix with
the entry $m,n$ denoted by $[\mathbb X(x)]_{n,m}$:
\begin{eqnarray}
\label{eq:DefBonD}
 B(x) = \bbmatrix{e^{-\mu_n x} b_{1n} & e^{\mu_n x} b_{2n}}, \quad \{b_{1n}\}, \{ b_{2n} \}\in\mathcal H, \\
\label{eq:DefXonD} 
[\mathbb X(x)]_{n,m} =  - \dfrac{b_{1n}\bar b_{1m} e^{-(\mu_n +\bar\mu_m) x} + b_{2n}\bar b_{2m} e^{(\mu_n +\bar\mu_m) x}}{2(\mu_n+\bar\mu_m)}.
\end{eqnarray}
we also assume that $ [\mathbb X(x)]_{n,m} = \dfrac{b_{1n}\bar b_{1m} - b_{2n}\bar b_{2m}}{2} x$, whenever $\mu_n +\bar\mu_m = 0$.
The fact that $B(x)$ is a well defined operator is immediate from the definition. Indeed, since the sequence $\mu_n$ is bounded the term
$e^{\pm \mu_n x}$ is uniformly bounded in absolute value by $e^{M x}$. The fact that the operator $\mathbb X(x)$ is bounded also easily follows
from the definitions and from the assumption on $D$:
\[ \| \mathbb X(x)\| \leq \dfrac{e^{2Mx}}{2} \| [\dfrac{|b_{1n}\bar b_{1m}| + |b_{2n}\bar b_{2m}|}{\mu_n+\bar\mu_m}] \| \leq
\dfrac{e^{2Mx}}{2} (\| b_{1n}\|^2  + \| b_{2n} \|^2) < \infty.
\]
Under condition that $\mathbb X_0 = \mathbb X(0)$ is invertible, there exists a non trivial interval (of length at least $\dfrac{1}{\|\mathbb X_0^{-1}\|}$) on which $\mathbb X(x)$
is invertible too. Thus we obtain the following Theorem.
\begin{thm} Suppose that we are given a bounded, separated from zero set $\{\mu_n\}$ and two $\ell^2$ sequences
$\{b_{1n}\}, \{ b_{2n} \}$. Define a collection 
\[ 
\mathfrak{V}_{NLS} = (A, B(x), \mathbb X(x); \sigma_1, 
\sigma_2, \gamma, \gamma_*(x);
\mathcal{H},\mathbb C^2),
\]
where $\mathcal H = \ell^2$, $A = \diag(2\mu_n)$, $\mathbb X(x)$ and $B(x)$ are defined by \eqref{eq:DefXonD} and by \eqref{eq:DefBonD} for
$|b_{1n}| + |b_{2n}| \neq 0$. Suppose that the operator $\mathbb X(0)$ is invertible.
Then the collection $\mathfrak{V}_{NLS}$ is an NLS regular vessel existing on a non trivial interval $\mathrm I$ including zero of length at least $\dfrac{1}{\|\mathbb X_0^{-1}\|}$.
\end{thm}

\section{Evolutionary regular NLS vessel}
We present a construction of solutions of the equation \eqref{eq:ENLS} which has initial value $\beta(x,0)$ arising from a regular
NLS vessel. For this we will insert dependence on the variable $t$ into the vessel operators and postulate evolution of the operators $B, \mathbb X$
with respect to $t$. This is done in the next Definition.
\begin{defn} The collection
\begin{equation} \label{eq:NLSEvolv}
\mathfrak{V}_{ENLS} = (A, B(x,t), \mathbb X(x,t); \sigma_1, 
\sigma_2, \gamma, \gamma_*(x);
\mathcal{H},\mathbb C^2), 
\end{equation}
is called \textbf{evolutionary regular NLS vessel} if additionally to equations \eqref{eq:DB}, \eqref{eq:Lyapunov}, \eqref{eq:DX}, \eqref{eq:Linkage}
the bounded operators satisfy also the equations
\begin{eqnarray}
\label{eq:DBt}  \dfrac{\partial}{\partial t} B(x,t) = i A \dfrac{\partial}{\partial x} B(x,t), \\
\label{eq:DXt} \dfrac{\partial}{\partial t} \mathbb X =  i A B \sigma_2 B^* - i B\sigma_2 B^* A^*.
\end{eqnarray}
\end{defn}
\begin{thm} Moments $H_N(x,t)$ of the vessel $\mathfrak{V}_{ENLS}$ satisfy the following recurrence equation
\[ (H_N)'_t = i (H_{N+1})'_x + i (H_0)'_x \sigma_1 H_N.
\]
and the transfer functions satisfies
\begin{equation}\label{eq:DSt} S'_t = i\lambda S'_x + i(H_0)'_x \sigma_1 S.
\end{equation}
\end{thm}
\noindent\textbf{Proof:} We will calculate the derivative of $H_N$ using the vessel conditions:
\[ \begin{array}{lllllllll}
(H_N)'_t & = (B^*\mathbb X^{-1}A^NB)'_t = (B^*)'_t\mathbb X^{-1}A^NB + B^*(\mathbb X^{-1})'_tA^NB + B^*\mathbb X^{-1}A^N(B)'_t \\
& = \text{using \eqref{eq:DBt}, \eqref{eq:DXt}} \\
& = (B^*)'_x(-iA^*)\mathbb X^{-1}A^NB - B^*\mathbb X^{-1}( i A B \sigma_2 B^* - i B\sigma_2 B^* A^*)\mathbb X^{-1}A^NB + B^*\mathbb X^{-1}A^N(iA)B'_x\\
& = \text{using \eqref{eq:Lyapunov}}\\
& = i (B^*)'_x(\mathbb X^{-1}A+\mathbb X^{-1}B\sigma_1B^*\mathbb X^{-1})A^NB - \\
 & \quad - i B^*\mathbb X^{-1} A B \sigma_2 B^* \mathbb X^{-1}A^NB - i B^*\mathbb X^{-1} B\sigma_2 B^* (\mathbb X^{-1}A+\mathbb X^{-1}B\sigma_1B^*\mathbb X^{-1}) A^NB + \\
 &\quad\quad + i B^*\mathbb X^{-1}A^{N+1}B'_x\\
 & = i [(B^*)'_x\mathbb X^{-1}A^{N+1}B -  B^*\mathbb X^{-1} A B \sigma_2 B^* \mathbb X^{-1}A^{N+1}B +i B^*\mathbb X^{-1}A^{N+1}B'_x] + \\
 & \quad + i (B^*)'_x\mathbb X^{-1}B\sigma_1B^*\mathbb X^{-1}A^NB -i B^*\mathbb X^{-1} A B \sigma_2 B^* \mathbb X^{-1}A^NB -
 i  B^*\mathbb X^{-1} B\sigma_2 B^* \mathbb X^{-1}B\sigma_1B^*\mathbb X^{-1} A^NB \\
 & = \text{using \eqref{eq:DX} and the definition of moments}\\
 & = i (H_{N+1})'_x + i (B^*)'_x\mathbb X^{-1}B\sigma_1 H_N + B^*\mathbb X^{-1}(B)'_x\sigma_1 H_N +  B^*(\mathbb X^{-1})'_x B\sigma_1H_N \\
 & = i (H_{N+1})'_x + i(H_0)'_x \sigma_1 H_N.
\end{array} \]
The formula for $S'_t$ is an immediate consequence of the series representation of $S$ using moments and equating the powers of $\lambda$.
\qed
\begin{thm} \label{thm:EbeatEqn}
Suppose that $\mathfrak{V}_{ENLS}$ is an evolutionary NLS vessel, then $\beta(x,t)$ \eqref{eq:DefBeta} satisfies the evolutionary NLS equation
\eqref{eq:ENLS}.
\end{thm}
\noindent\textbf{Proof:} Consider the equality of mixed derivatives
\[ \dfrac{\partial}{\partial x}\dfrac{\partial}{\partial t} S = \dfrac{\partial}{\partial t}\dfrac{\partial}{\partial x} S.
\]
Plugging \eqref{eq:DS} and \eqref{eq:DSt} here and representing $S$ as a series of moments, we can further consider
both sides of this equality and make the corresponding coefficients of $\lambda^{-n}$ equal:
\[ \dfrac{\partial}{\partial x}  [i\lambda S'_x + i(H_0)'_x \sigma_1 S] = 
\dfrac{\partial}{\partial t} [\sigma_1^{-1}(\sigma_2 \lambda + \gamma_*) S - S \sigma_1^{-1}(\sigma_2 \lambda + \gamma)].
\]
For example, taking the free coefficient, we wil obtain that
\begin{equation} \label{eq:Dgamma*tKdV}
(\gamma_*)'_t = - i \gamma_* (H_0)'_x\sigma_1 + i \sigma_1 (H_0)''_{xx} \sigma_1 +i \sigma_1 (H_0)'_x \gamma_*.
\end{equation} 
Finally, From the linkage condition \eqref{eq:Linkage}
\[
\gamma_* = \bbmatrix{0&\beta(x,t)\\-\beta^*(x,t)&0}
\]
and using the formula for $(H_0)'_x = \bbmatrix{|\beta|^2 & (\beta)'_x \\(\beta^*)'_x&-|\beta|^2}$
we will obtain that the 12 entry is translated into equation \eqref{eq:ENLS} for $\beta(x,t)$, defined by \eqref{eq:DefBeta}.\qed

So in order to solve \eqref{eq:ENLS} with initial $\beta(x,0)$ which arises from a vessel, it is enough to add dependence on $t$ so that
$B(x,0) = B(x)$, $\mathbb X(x,0) = \mathbb X(x)$ and the differential equations \eqref{eq:DBt}, \eqref{eq:DXt} hold. We will show that it is 
a simple task in the next examples.

\section{Examples of constructions of solutions of the evolutionary NLS equation}
We present examples of solutions of the evolutionary NLS \eqref{eq:ENLS} when the initial value for $t=0$,
$\beta(x)$ is analytic on $\mathbb R$ and arises from a regular NLS vessel. We show how to construct the evolutionary
vessel, coinciding with the vessel realizing $\beta(x)$ for $t=0$. Then the beta function of this evolutionary vessel is a solution
of \eqref{eq:ENLS} with the initial value.

\subsection{\label{sec:EVesselFromS0}Construction of a solution for evolutionary NLS vessel from a realized function}
Suppose that $\beta(x)$ was constructed from a realized function as in Section \ref{sec:VesselFromS0}
\[ S(\lambda,x_0) = I - B_0^* \mathbb X_0^{-1}(\lambda I - A)^{-1} B_0.
\]
Then the construction can proceed the following steps, each one requiring a solution of Linear Differential
equation with initial value. Construct $B(x)$ and $\mathbb X(x)$ by formulas \eqref{eq:DefBFromS0}, \eqref{eq:DefXFromS0}. Then solve
\[ \dfrac{\partial}{\partial t} B(x,t) = i A \dfrac{\partial}{\partial x} B(x,t) = i A (-A B(x,t) \sigma_2), B(x,0) = B(x).
\]
Finally, define
\[ \mathbb X(x,t) = \mathbb X(x) + \int_0^t [i A B(x,s) \sigma_2 B^*(x,s) - i B(x,s) \sigma_2 B^*(x,s) A^*] ds.
\]
All the vessel equations will be satisfied by the construction and can be easily verified.

\subsection{\label{sec:EVesselFromGamma} Solution of the evolutionary NLS with the spectrum on a curve $\Gamma$}
If we are given $\beta(x)$ which arises from a vessel on a curve, we can solve explicitly equations for $B(x,t)$ and
$\mathbb X(x,t)$ as follows:
\begin{equation} \label{eq:DefEBonGamma}
 B(x) = \bbmatrix{\exp(-\mu x -2i\mu^2 t) b_1(\mu) & \exp(\mu x + 2i\mu^2 t) b_2(\mu)}, \quad b_1(\mu), b_2(\mu)\in\mathcal H.
\end{equation}
This function coincides with $B(x)$ for $t=0$ and satisfies \eqref{eq:DB} and \eqref{eq:DBt}. The formula for $\mathbb X(x,t)$ is as follows:
\begin{multline} \label{eq:DefEXonGamma}
(\mathbb X(x) f) (\mu) = - \dfrac{1}{2} \int_a^b \dfrac{b_1(\mu) \bar b_1(\mu(s)) \exp(-(\mu+\bar\mu(s))x - i (\mu^2+\bar\mu^2(s))t)}{\mu + \bar\mu(s)} f(\mu(s)) ds
- \\
- \dfrac{1}{2} \int_a^b \dfrac{ b_2(\mu) \bar b_2(\mu(s)) \exp((\mu+\bar\mu(s))x + i(\mu^2 + \mu^2(s))t)}{\mu + \bar\mu(s)} f(\mu(s)) ds. 
\end{multline}
This constructs a solution of the \eqref{eq:ENLS} with initial $\beta(x)$ arising from an NLS vessel with the spectrum on a curve $\Gamma$

\subsection{\label{sec:EVesselFromD}Solution of the evolutionary NLS with the spectrum on a discrete set}
Similarly to the continuous spectrum case we define
\begin{eqnarray}
\label{eq:DefEBonD}
 B(x) = \bbmatrix{\exp(-\mu_n x-2i\mu^2_nt) b_{1n} & \exp(\mu_n x + 2i\mu^2_n) b_{2n}}, \quad \{b_{1n}\}, \{ b_{2n} \}\in\mathcal H, \\
\label{eq:DefEXonD} 
[\mathbb X(x)]_{n,m} =  - \dfrac{b_{1n}\bar b_{1m} e^{-(\mu_n +\bar\mu_m) x - 2i(\mu_n^2 +\bar\mu_m^2)t} + b_{2n}\bar b_{2m} e^{(\mu_n +\bar\mu_m) x + 2i(\mu_n^2 +\bar\mu_m^2)t}}{2(\mu_n+\bar\mu_m)}.
\end{eqnarray}
we also assume that $ [\mathbb X(x)]_{n,m} = \dfrac{b_{1n}\bar b_{1m} - b_{2n}\bar b_{2m}}{2} (x + 2 i \mu_n t)$, whenever $\mu_n +\bar\mu_m = 0$.

\subsection{\label{sec:Solitons}Solitons}
We can also consider finite dimensional case: $\dim\mathcal H<\infty$. Suppose that $\mathcal H = \mathbb C^N$ and fixing non zero values
$\mu_1,\mu_2,\ldots,\mu_N$, satisfying $\mu_i + \bar\mu_j\neq 0$ for all $1\leq i,j\leq N$, we will obtain that
\[ B(x) = \bbmatrix{
\exp(-\mu_1 x-2i\mu^2_1t) b_{11} & \exp(\mu_1 x + 2i\mu^2_1) b_{21} \\
\exp(-\mu_2 x-2i\mu^2_2t) b_{12} & \exp(\mu_1 x + 2i\mu^2_2) b_{22} \\
\vdots & \vdots \\
\exp(-\mu_N x-2i\mu^2_Nt) b_{1N} & \exp(\mu_N x + 2i\mu^2_N) b_{2N} 
}
\]
and the operator $\mathbb X(x)$ is $N\times N$ matrix, defined by the formula \eqref{eq:DefEXonD}. Then the function
\[
\beta(x,t) = \bbmatrix{1&0} B^*(x) \mathbb X^{-1}(x) B(x)\bbmatrix{0\\1}
\]
is a Soliton solution of \eqref{eq:ENLS}, because it is constructed from pure exponents.

\bibliographystyle{alpha}
\bibliography{../../biblio}

\end{document}